\pdfoutput=1
\documentclass{article}
\usepackage[a4paper, margin=1.5in]{geometry}

\usepackage[T1]{fontenc}      % Encoding
\usepackage[utf8]{inputenc}   % Encoding
\usepackage{microtype}        % Improved typesetting
\usepackage{mlmodern}

\usepackage{amsmath}
\usepackage{amsthm}
\usepackage{amssymb}
\usepackage{mathtools}
\usepackage{thmtools}
\usepackage{commath}
\usepackage{bbm}
\usepackage{xfrac}
\usepackage[numbers]{natbib}
\usepackage{cleveref}

\usepackage{xcolor}

\usepackage{tikz}
\usepackage{pgfplots}
\usepgfplotslibrary{groupplots}
\pgfplotsset{compat=1.18}

\numberwithin{equation}{section}
\declaretheorem[Refname={Theorem,Theorems}]{theorem}
\numberwithin{theorem}{section} % This should commented out if there is to be a non-section-based nvumbering

\declaretheorem[name=Proposition,numberlike=theorem,Refname={Proposition,Propositions}]{proposition}

\newcommand{\R}{\mathbb{R}} % Set of real numbers
\newcommand{\N}{\mathbb{N}} % Set of natural numbers
\newcommand{\C}{\mathbb{C}} % Set of complex numbers

\DeclarePairedDelimiterX\Set[2]{\lbrace}{\rbrace}%
{ #1 \,:\, #2 }

\newcommand{\err}{\mathrm{err}} % Set of complex numbers
\DeclareMathOperator{\supp}{supp}

\renewcommand{\i}{\mathrm{i}}

\begin{document}

\title{\text{Bumps and polynomials in RKHSs} \\ \text{of translation-invariant kernels}}

\author{Toni Karvonen \vspace{0.2cm} \\ \emph{School of Engineering Sciences} \\ \emph{Lappeenranta--Lahti University of Technology LUT, Finland.}}

\maketitle

\begin{abstract}
  \noindent We use the uncertainty principle of harmonic analysis and the Beurling--Malliavin theorem to describe when a reproducing kernel Hilbert space of a translation-invariant positive-definite kernel contains bump functions and polynomials.
  How fast the spectral density of the kernel decays controls whether or not these functions are included in the Hilbert space.
\end{abstract}

\section{Introduction}

In the broadest sense of the term, a bump function is any non-trivial compactly supported function (see \Cref{fig:bumps}). 
\Cref{thm:wce} shows that bump functions play a crucial role in determining whether or not approximating functions in a \emph{reproducing kernel Hilbert space} (RKHS) at a point $x$ requires information at points close to $x$.
Recall that, by the Moore--Aronszajn theorem~\cite[p.\@~19]{Berlinet2004}, every positive-semidefinite kernel $K \colon X \times X \to \R$ on a set $X$ induces a unique RKHS, $H(X)$, of functions $f \colon X \to \R$ in which $K$ is reproducing, meaning that $\langle f, K(\cdot, x) \rangle_{H(X)} = f(x)$ for every $f \in H(X)$ and $x \in X$.
The \emph{worst-case error}, a staple of information-based complexity~\cite{NovakWozniakowski2008},
\begin{equation*}
  \err_n(x) = \sup_{ 0 \neq f \in H(X) } \inf_{u_1, \ldots, u_n \in \R} \frac{ \lvert f(x) - \sum_{i=1}^n u_i f(x_i) \rvert  }{\lVert f \rVert_{H(X)}}
\end{equation*}
quantifies how easy it is to approximate functions in $H(X)$ at a point $x \in X$ using linear combinations of function evaluations at points $x_1, \ldots, x_n$.
The worst-case error is known as power function in scattered data approximation~\cite{Wendland2005} and as posterior standard deviation in Gaussian process modelling literature~\cite{RasmussenWilliams2006}.
Throughout this article $B(x, \varepsilon)$ denotes the closed ball of radius $\varepsilon$ centered at $x$.

\pgfplotsset{colormap={graybg}{rgb(0cm)=(0.6,0.6,0.6); rgb(1cm)=(0,0,0)}}
\begin{figure}
\centering
\begin{tikzpicture}
\begin{groupplot}[
    group style={group size=3 by 1, horizontal sep=0.5cm},
    width=6cm, height=6cm,
]

% --- 1D bump function ---
\nextgroupplot[
    hide axis,
    xmin=-2, xmax=2,
    ymin=-0.2, ymax=1.3,
    samples=200,
    domain=-2:2,
    smooth,
]
\addplot[ultra thick, black] {(abs(x)<1) * exp(-1/(1-x^2))*2.71828};

% --- 1D bump function ---
\nextgroupplot[
    hide axis,
    xmin=-2, xmax=2,
    ymin=-0.2, ymax=1.3,
    samples=200,
    domain=-2:2,
    smooth,
]
\addplot[ultra thick, black] {(abs(x)<1) * exp(-1/(1-x^4))*2.71828};

% --- 2D bump function ---
\nextgroupplot[
    hide axis,
    domain=-1.5:1.5,
    y domain=-1.5:1.5,
    samples=40,
    samples y=40,
    axis equal image,
    colormap name=graybg,
]
\addplot3[mesh, thin, shader=interp] {(x^2+y^2<1) * exp(-1/(1-(x^2+y^2)^4))*2.71828};

\end{groupplot}
\end{tikzpicture}
\caption{The bump function $\psi(x) = \exp(-1/(1-\lVert x \rVert^{2\alpha})) \mathbbm{1}_{B(0, 1)}(x)$ supported on the closed unit ball $B(0, 1)$ on $\R$ with $\alpha = 1$ and $\alpha = 2$ and on $\R^2$ with $\alpha = 4$.}
\label{fig:bumps}
\end{figure}

\begin{theorem} \label{thm:wce}
  Let $X$ be a metric space.
  Suppose that $K(\cdot, x)$ is continuous for every $x \in X$ and that $x \mapsto K(x, x)$ is locally bounded.
  \begin{itemize}
    \item[(i)] If every closed ball $B(x_0, \varepsilon) \subseteq X$ there is a non-trivial function $\psi \in H(X)$ such that $\supp(\psi) \subseteq B(x_0, \varepsilon)$ and $\psi(x_0) \neq 0$, then $\err_n(x) \to 0$ as $n \to \infty$ if and only if $x \in \operatorname{cl} \{ x_i \}_{i=1}^\infty$.
    \item[(ii)] If no non-trivial function in $H(X)$ vanishes on an open set, then $\err_n(x) \to 0$ for any $x \in X$ whenever $\operatorname{cl}\{x_i\}_{i=1}^\infty$ has non-empty interior.
  \end{itemize}
\end{theorem}
\begin{proof}
  See \Cref{sec:proofs}.
\end{proof}

Sobolev spaces are no doubt the most familiar function spaces that contain bump functions (e.g., the functions in \Cref{fig:bumps}) and fall under Part~(i) of \Cref{thm:wce}, whereas no space that consists of analytic functions can contain bump functions by the uniqueness of analytic continuation.
Moreover, if $H(X)$ consists of analytic functions, in Part~(ii) the requirement that the closure have non-empty interior can be relaxed to the set having an accumulation point if $X = \R$ or the closure having a positive Lebesgue measure if $X = \R^d$ for $d \geq 2$; see~\cite[Corollary 1.2.7]{KrantzParks2002} and~\cite{Mityagin2020}.
It is worth noting that in optimisation literature~\cite{VazquezBect2010a} a positive-definite kernel is sometimes said to have the no-empty-ball property if $\err_n(x) \to 0$ as $n \to \infty$ is equivalent to $x \in \operatorname{cl} \{ x_i \}_{i=1}^\infty$.

The purpose of this article is to use results on the \emph{uncertainty principle} of harmonic analysis~\cite{HavinJoricke1994} to describe when RKHSs of translation-invariant kernels contain bump functions and when they do not in terms of the rate of decay of the spectral density of the kernel.
As a by-product we describe RKHSs of translation-invariant kernels that contain polynomials, a problem of some independent interest~\cite{Minh2010, Steinwart2006}.
The uncertainty principle states that a non-trivial function and its Fourier transform cannot be simultaneously ``too small''.
Let $f \colon \R \to \R$ have Fourier transform $\mathcal{F} f$.
Smallness is controlled by the \emph{logarithmic integral}
\begin{equation*}
  \mathcal{L}(\mathcal{F} f) = \int_\R \frac{\log \lvert \mathcal{F} f(\xi) \rvert }{1 + \xi^2} \dif \xi .
\end{equation*}
A function $f \colon \R \to \R$ has \emph{semibounded support} if there is $a \in \R$ such that either $\supp(f) \subseteq (-\infty, a]$ or $\supp(f) \subseteq [a, \infty)$.
A classical result in harmonic analysis (see \Cref{thm:semibounded}) states that a function $f$ with semibounded support is identically zero if $\mathcal{L}(\mathcal{F} f) = -\infty$.
This happens, for instance, if $\mathcal{F} f(\xi) = \exp(-\lvert \xi \rvert^\alpha)$ for $\alpha \geq 1$.
That is, if its Fourier transform decays too fast, a function cannot be compactly supported unless it is zero.
The Beurling--Malliavin theorem (see \Cref{thm:BM}) provides a converse to the above result.
Namely, if $\eta$ is a non-negative Lipschitz function such that $\mathcal{L}(\eta) > -\infty$, then there exists a non-trivial bump function $f$ such that $\mathcal{F} f(\xi) = O(\eta(\xi))$.
Our contribution consists in applying these theorems of harmonic analysis to translation-invariant kernels.

\section{Results}

Let $\Phi \in C(\R^d) \cap L^1(\R^d)$ be a real-valued positive-definite function on $\R^d$ and $K(x, y) = \Phi(x - y)$ the corresponding translation-invariant positive-definite kernel.
It is a well-known result~\cite[Theorem~10.12]{Wendland2005} that goes back at least to the work of Kimeldorf and Wahba~\cite{KimeldorfWahba1970} that the RKHS of such a translation-invariant kernel is characterised in terms of the Fourier transform, or \emph{spectral density},
\begin{equation*}
  \mathcal{F}\Phi(\xi) = \frac{1}{(2\pi)^{d/2}} \int_{\R^d} \Phi(x) e^{-\i x \cdot \xi} \dif x
\end{equation*}
of the kernel.
That is,
\begin{equation} \label{eq:rkhs}
    H(\R^d) = \Set[\bigg]{ f \in C(\R^d) \cap L^2(\R^d) }{ \int_{\R^d} \frac{\lvert \mathcal{F}f(\xi) \rvert^2}{\mathcal{F}\Phi(\xi)} \dif \xi < \infty } .
\end{equation}
The following theorem states that $H(\R)$ does not contain compactly supported functions if the spectral density decays too fast.

\begin{theorem}[Uncertainty principle for kernels] \label{thm:main-theorem-converse}
  Let $d = 1$.
  Suppose that $\mathcal{F} \Phi(\xi) = O( e^{-\omega(\xi)} )$ for a function $\omega \colon \R \to [0, \infty)$ such that
  \begin{equation} \label{eq:ass-3}
    \int_\R \frac{\omega(\xi)}{1 + \xi^2} \dif \xi = \infty .
  \end{equation}
  If $f \in H(\R)$ has semibounded support, then $f \equiv 0$.
\end{theorem}
\begin{proof}
  See \Cref{sec:proofs}.
\end{proof}

A somewhat less general version is available for arbitrary $d$.

\begin{theorem}[Paley--Wiener for kernels] \label{thm:analytic}
  Let $d \geq 1$. If $\mathcal{F}\Phi(\xi) = O(e^{-\alpha \lVert \xi \rVert})$ for some $\alpha > 0$, then every function in $H(\R^d)$ is real analytic and the zero function is the only function in $H(\R^d)$ that vanishes on an open set.
\end{theorem}
\begin{proof}
  See \Cref{sec:proofs}.
\end{proof}

The following theorem describes when bump functions are contained in the RKHS of a translation-invariant kernel.
The theorem is available in the multivariate setting due to Vasilyev and Bergman \cite{Bergman2025,Vasileyv2025} who have recently proved a multivariate generalisation (see \Cref{thm:BM}) of the Beurling--Malliavin theorem from 1962~\cite{BeurlingMalliavin1962}.
A function $\omega$ defined on a normed space is \emph{radial} if $\omega(x)$ depends only on $\lVert x \rVert$. 

\begin{theorem}[Beurling--Malliavin for kernels] \label{thm:main-theorem}
  Let $d \geq 1$.
  If there is a radial Lipschitz function $\omega \colon \R^d \to [0, \infty)$ and a constant $c > 0$ such that
  \begin{equation} \label{eq:ass-1}
    e^{-\omega(\xi)} \leq c \cdot \mathcal{F}\Phi(\xi) 
  \end{equation}
  for all $\xi \in \R^d$ and 
  \begin{equation} \label{eq:ass-2}
    \int_{\R^d} \frac{\omega(\xi)}{(1 + \lVert \xi \rVert^2 )^{(d+1)/2}} \dif \xi < \infty,
  \end{equation}
  then for every closed ball $B(x, \varepsilon)$ there is $\psi \in H(\R^d)$ such that $\supp(\psi) \subseteq B(x, \varepsilon)$, $\psi(x) = 1$, and $\int_{\R^d} \psi(x) \dif x \neq 0$.
\end{theorem}
\begin{proof}
  See \Cref{sec:proofs}.
\end{proof}

The integrability conditions~\eqref{eq:ass-3} and~\eqref{eq:ass-2} also (essentially) control whether or not an RKHS on a bounded subset contains polynomials.
The RKHS of $K$ on an arbitrary set $\Omega \subseteq \R^d$ is the restriction
\begin{equation*}
  H(\Omega) = \Set{ f|_\Omega }{ f \in H(\R^d)} .
\end{equation*}

\begin{theorem}[Polynomials in an RKHS] \label{thm:polynomials}
  Let $d \geq 1$ and $\Omega \subseteq \R^d$.
  \begin{enumerate}
    \item[(i)] Suppose that $\Omega$ is unbounded. Then $f(x) \to 0$ as $\lVert x \rVert \to \infty$ if $f \in H(\Omega)$. In particular, $H(\Omega)$ contains no non-trivial polynomials.
    \item[(ii)] Suppose that $\Omega$ is bounded and has a non-empty interior. If $\mathcal{F}\Phi(\xi) = O(e^{-\alpha \lVert \xi \rVert})$ for some $\alpha > 0$, then $H(\Omega)$ contains no non-trivial polynomials.
    \item[(iii)] Suppose that $\Omega$ is bounded. If there is a radial Lipschitz function $\omega \colon \R^d \to [0, \infty)$ and a constant $c > 0$ such that~\eqref{eq:ass-1} and~\eqref{eq:ass-2} hold, then $H(\Omega)$ contains all polynomials.
  \end{enumerate}  
\end{theorem}
\begin{proof}
  See \Cref{sec:proofs}.
\end{proof}

We note that Theorem~3 in~\cite{SunZhou2008} closely resembles Part~(ii) of \Cref{thm:polynomials}.

\section{Examples}

This section contains a collection of examples that illustrate \Cref{thm:analytic,thm:main-theorem,thm:main-theorem-converse}.

\subsection{Matérn kernels}

Let $\nu > 0$ be a regularity and $\ell > 0$ a scale parameter.
The \emph{Matérn kernel} of order $\nu$ is defined by the function
\begin{equation*}
  \Phi(z) = \frac{2^{1-\nu}}{\Gamma(\nu)} \bigg( \frac{\sqrt{2\nu} \lVert z \rVert }{\ell} \bigg)^\nu \mathcal{K}_\nu \bigg( \frac{\sqrt{2\nu} \lVert z \rVert }{\ell} \bigg),
\end{equation*}
where $\mathcal{K}_\nu$ is the modified Bessel function of the second kind of order $\nu$.
Matérns are particularly popular in Gaussian process modelling.
The spectral density decays polynomially~\citep[Theorem~6.13]{Wendland2005}:
\begin{equation*}
  \mathcal{F} \Phi (\xi) = \frac{2^{d/2} \Gamma(\nu+d/2)}{\Gamma(\nu)} \bigg( \frac{2\nu}{\ell^2} \bigg)^\nu
\bigg(\frac{2\nu}{\ell^2}+ \lVert \xi \rVert^2\bigg)^{-(\nu+d/2)} .
\end{equation*}
From~\eqref{eq:rkhs} it is then immediate that $H(\R^d)$ is a Bessel potential space of order $\nu + d/2$ and norm-equivalent to the standard Sobolev space $W^{\alpha,2}(\R^d)$ when $\alpha = \nu + d/2$ is an integer.
It is a basic result that Sobolev spaces contain bump functions, standard examples being the infinitely differentiable bump function $\psi(x) = \exp(-1/(1-\lVert x \rVert^2)) \mathbbm{1}_{B(0,1)}(x)$ and its translations.
One could also use \Cref{thm:main-theorem} to verify this by setting
\begin{equation*}
  \omega(\xi) = (\nu + d/2) \log\bigg( \frac{2\nu}{\ell^2} + \lVert \xi \rVert^2 \bigg),
\end{equation*}
which is Lipschitz and clearly satisfies~\eqref{eq:ass-2}.

\subsection{Gevrey kernels}

The remaining examples have power exponential spectral densities.
Let $\alpha > 0$.
A kernel with spectral density
\begin{equation} \label{eq:exp-kernel}
  \mathcal{F} \Phi(\xi) = \exp(- \alpha \lVert \xi \rVert^\beta)
\end{equation}
is a \emph{Gevrey kernel} if $\beta \in (0, 1)$.
It is not difficult to show that these kernels induce RKHSs that consist of functions that lie in the Gevrey class of order $\sigma = 1/\beta$. 
Gevrey classes contain functions that are infinitely differentiable but need not be analytic~\cite[Chapter~I]{Rodino1993}.
The function
\begin{equation*}
  \omega(\xi) = \alpha ( 1 + \lVert \xi \rVert)^\beta 
\end{equation*}
is Lipschitz and satisfies~\eqref{eq:ass-1} and~\eqref{eq:ass-2}.
Therefore $H(\R^d)$ contains bump functions by \Cref{thm:main-theorem}.
To find the inverse Fourier transform of~\eqref{eq:exp-kernel} is non-trivial when $\beta \in (0, 1)$. 
The kernel is available in certain special cases because these inverse Fourier transforms correspond to density functions of stable distributions~\cite{GaroniFrankel2002, Zolotarev1986}.
Let $d = 1$.
The simplest of these special cases is $\beta = 2/3$.
The resulting kernel, which we call the \emph{Whittaker kernel}, is 
\begin{equation*}
  \Phi(z) = \frac{1}{\sqrt{6} \, \lvert z \rvert} \exp\bigg(\frac{2}{27 \gamma^2 z^2} \bigg) \mathcal{W}_{-1/2, 1/6} \bigg(\frac{4}{27 \gamma^2 z^2} \bigg),
\end{equation*}
where $\gamma = \alpha^{-1/\beta}$ and $\mathcal{W}_{\kappa, \mu}$ is the Whittaker function, which has the expression
\begin{equation*}
  \mathrm{W}_{\kappa, \mu}(z) = \frac{z^\kappa e^{-z/2}}{\Gamma(1/2 - \kappa + \mu)} \int_0^\infty t^{-\kappa-1/2+\mu} \bigg(1 + \frac{t}{z} \bigg)^{\kappa-1/2+\mu} e^{-t} \dif t
\end{equation*}
when $\kappa-1/2-\mu \in \C$ is not an integer and has non-positive real part.

\subsection{Inverse quadratic kernel}

Setting $\beta = 1$ in~\eqref{eq:exp-kernel} yields $\mathcal{F}\Phi(\xi) = \exp(-\alpha \lVert \xi \rVert)$.
When $d = 1$ the inverse Fourier transform is simply
\begin{equation*}
  \Phi(z) = \sqrt{\frac{2}{\pi}} \cdot \frac{\alpha}{\alpha^2 + z^2} .
\end{equation*}
This defines the \emph{inverse quadratic kernel} popular in scattered data approximation.
Because 
\begin{equation*}
  \int_\R \frac{\lvert \xi \rvert^\beta}{1 + \xi^2} \dif \xi < \infty
\end{equation*}
if and only if $\beta \in (0, 1)$, the inverse quadratic kernel is a boundary case among kernels whose spectral densities have the form~\eqref{eq:exp-kernel}.
By \Cref{thm:main-theorem-converse,thm:analytic} its RKHS does not contain bump functions and consists of analytic functions.
 
\subsection{Gaussian kernel}

Setting $\beta = 2$ in~\eqref{eq:exp-kernel} yields the Gaussian spectral density $\mathcal{F}\Phi(\xi) = \exp(-\alpha \lVert \xi \rVert^2)$ and the \emph{Gaussian kernel}
\begin{equation*}
  \Phi(z) = (2\alpha)^{-d/2} \exp\bigg(\! - \frac{\lVert z \rVert^2}{4\alpha}\bigg) ,
\end{equation*}
which is particularly popular in machine learning.
By \Cref{thm:analytic} the RKHS of the Gaussian kernel does not contain bump functions and consists of analytic functions.
The same applies to all spectral densities of the form~\eqref{eq:exp-kernel} with $\beta \geq 1$.

\section{Proofs} \label{sec:proofs}

\begin{proof}[Proof of \Cref{thm:wce}]
  Consider Part~(i). 
  Assume that $x \in \operatorname{cl}\{x_i\}_{i=1}^\infty$.
  For any $m \leq n$ we have
  \begin{equation*}
    \begin{split}
    \err_{n}(x) &= \sup_{ 0 \neq f \in H(X) } \inf_{u_1, \ldots, u_n \in \R} \frac{ \lvert f(x) - \sum_{i=1}^n u_i f(x_i) \rvert  }{\lVert f \rVert_{H(X)}} \\
    &\leq \sup_{ 0 \neq f \in H(X) } \inf_{u_m \in \R} \frac{ \lvert f(x) - u_m f(x_m) \rvert }{\lVert f \rVert_{H(X)}} \\
    &\leq \sup_{ 0 \neq f \in H(X) } \inf_{u_m \in \R} \frac{ \lvert \langle f, K(\cdot, x) - u_m K(\cdot, x_m) \rangle_{H(X)} \rvert }{\lVert f \rVert_{H(X)}} \\
    &\leq \inf_{u_m \in \R} \sqrt{K(x, x) - 2 u_m K(x, x_m) + u_m^2 K(x_m, x_m)} \\
    &= \sqrt{K(x, x) - K(x, x_m)^2 / K(x_m, x_m)} ,
    \end{split}
  \end{equation*}
  where we used the reproducing property and the Cauchy--Schwarz inequality.
  Because $m \leq n$ is arbitrary and $K(\cdot, x)$ is continuous for every $x \in X$, $\err_n(x) \to 0$ follows from the assumption $x \in \operatorname{cl}\{x_i\}_{i=1}^\infty$.
  Assume then that $x \notin \operatorname{cl}\{x_i\}_{i=1}^\infty$.
  Then there is $\varepsilon > 0$ such that $B(x, \varepsilon)$ contains no $x_i$.
  Let $\psi \in H(X)$ be a non-trivial function such that $\supp(\psi) \subseteq B(x, \varepsilon)$ and $\psi(x) \neq 0$.
  By the definition of the worst-case error,
  \begin{equation*}
    \err_n(x) \geq \inf_{u_1, \ldots, u_n \in \R} \frac{ \lvert \psi(x) - \sum_{i=1}^n u_i \psi(x_i) \rvert  }{\lVert \psi \rVert_{H(X)}} = \frac{\lvert \psi(x) \rvert}{\lVert \psi \rVert_{H(X)}} > 0
  \end{equation*}
  for every $n$ because $\psi(x_i) = 0$ for every $i$.
  Therefore $\err_n(x)$ does not tend to zero.
  This proves the first claim.

  Consider Part~(ii).
  By assumption there is a non-empty closed ball $B$ such that $B \subseteq \operatorname{cl}\{x_i\}_{i=1}^\infty$.
  The assumptions on $K$ imply that all functions in $H(X)$ are continuous~\cite[p.\@~34]{Berlinet2004}.
  Because $\{x_i\}_{i=1}^\infty$ is dense in $B$, it follows from continuity that $f|_B \equiv 0$ for every $f$ in the orthogonal complement of $V = \operatorname{span}\{K(\cdot, x_i)\}_{i=1}^\infty$. 
  Since the only function in $H(X)$ to vanish on a non-empty ball is the zero function, we conclude that the orthogonal complement of $V$ is trivial.
  Therefore $V$ is dense in $H(X)$.
  By the reproducing property and the Cauchy--Schwarz inequality,
  \begin{equation*}
    \err_n(x) \leq \inf_{u_1, \ldots, u_n \in \R} \bigg\lVert K(\cdot, x) - \sum_{i=1}^n u_i K(\cdot, x_i) \bigg\rVert_{H(X)} .
  \end{equation*}
  This concludes the proof because $V$ is dense in $H(X)$ and $K(\cdot, x) \in H(X)$. 
\end{proof}

\Cref{thm:analytic} is nothing but an application to translation-invariant kernels of the Paley--Wiener theorem, which states that a function with an exponentially decaying Fourier transform is analytic.
We use the following simplified version of the Paley--Wiener theorem from~\cite[Theorem~IX.13]{ReedSimon1975}.

\begin{theorem}[Paley--Wiener I]
  If $f \in L^2(\R^d)$ and $\int_{\R^d} e^{\alpha \lVert \xi \rVert} \lvert \mathcal{F}f(\xi) \rvert^2 \dif \xi < \infty$ for some $\alpha > 0$, then $f$ is real analytic.
\end{theorem}

\begin{proof}[Proof of \Cref{thm:analytic}]
  Let $f \in H(\R^d)$.
  The characterisation~\eqref{eq:rkhs} and the assumption on $\mathcal{F}\Phi$ yield
  \begin{equation*}
    \infty > \int_{\R^d} \frac{\lvert \mathcal{F}f(\xi) \rvert^2}{\mathcal{F}\Phi(\xi)} \dif \xi \geq c \int_{\R^d} e^{\alpha \lVert \xi \rVert} \lvert \mathcal{F}f(\xi) \rvert^2 \dif \xi
  \end{equation*}
  for some $c > 0$.
  The Paley--Wiener theorem then implies that $f$ is real analytic.
  Finally, the zero function is the only real analytic function that vanishes on an open set~\cite{Mityagin2020}.
\end{proof}

The proof of \Cref{thm:main-theorem-converse} is based on the following classical theorem that can be found, for example, on page~36 of~\cite{HavinJoricke1994} (see also Chapter~2 of Part Two).
This theorem is one of the simplest manifestations of the uncertainty principle of harmonic analysis.

\begin{theorem} \label{thm:semibounded}
  If $f \in L^2(\R)$ has semibounded support and 
  \begin{equation*}
    \int_\R \frac{\log\lvert \mathcal{F} f(\xi) \rvert}{1 + \xi^2} \dif \xi = -\infty,
  \end{equation*}
  then $f \equiv 0$.
\end{theorem}

\begin{proof}[Proof of \Cref{thm:main-theorem-converse}]
    Let $f \in H(\R)$.
   The non-negative function
   \begin{equation*}
      g(\xi) = e^{\omega(\xi)} \lvert \mathcal{F}f(\xi) \rvert^2
   \end{equation*}
   is integrable by the characterisation~\eqref{eq:rkhs} and the assumption $\mathcal{F}\Phi(\xi) = O(e^{-\omega(\xi)})$. 
   Because $\log g \leq g$ and $g$ is integrable, $\int_\R \log g(\xi) / (1 + \xi^2) \dif \xi < \infty$.
   Therefore
   \begin{equation*}
      \int_\R \frac{\log \lvert \mathcal{F}f(\xi) \rvert }{1 + \xi^2} \dif \xi = -\frac{1}{2} \int_\R \frac{\omega(\xi) }{1 + \xi^2} \dif \xi + \frac{1}{2} \int_\R \frac{\log g(\xi) }{1 + \xi^2} \dif \xi = -\infty 
   \end{equation*}
   by~\eqref{eq:ass-3}.
   From \Cref{thm:semibounded} if follows that $f$ having semibounded support implies that $f \equiv 0$.
\end{proof}

\Cref{thm:main-theorem} is essentially a kernel version of the Beurling--Malliavin theorem~\cite{BeurlingMalliavin1962}.
The following multivariate version of the Beurling--Malliavin theorem was recently proved by Vasilyev and Bergman \cite{Bergman2025,Vasileyv2025}.
We also refer to~\cite{Mashreghi2006} for a non-specialist survey of this important theorem.

\begin{theorem}[Beurling--Malliavin for $d \geq 1$] \label{thm:BM}
  If $\omega \colon \R^d \to [0, \infty)$ is a radial Lipschitz function such that 
  \begin{equation*}
    \int_{\R^d} \frac{\omega(\xi)}{(1 + \lVert \xi \rVert^2 )^{(d+1)/2}} \dif \xi < \infty,
  \end{equation*}
  then for every $\varepsilon > 0$ there exists a non-trivial function $\psi \in L^2(\R^d)$ such that $\supp(\psi) \subseteq B(0, \varepsilon)$ and $\lvert \mathcal{F}\psi(\xi) \rvert \leq e^{-\omega(\xi)}$ for all $\xi \in \R^d$.
\end{theorem}

To ensure that the bump functions do not integrate to zero (a crucial ingredient in the proof of \Cref{thm:polynomials}) we use another version of the Paley--Wiener theorem.
This version states that a function is compactly supported if and only if its Fourier transform is an entire function of exponential type~\cite[Theorem~7.2.1]{Strichartz1994}.
Recall that a function $f \colon \C \to \C$ is \emph{entire} if can be expressed as a power series
\begin{equation*}
  f(z) = \sum_{n=0}^\infty a_n z^n
\end{equation*}
that converges on the whole of $\C$.
An entire function is of \emph{exponential type} $R > 0$ if 
\begin{equation*}
  f(z) = O(e^{R \lvert z \rvert}) .
\end{equation*}

\begin{theorem}[Paley--Wiener II]
  Let $R > 0$.
  Then $\mathcal{F}f$ is an entire function of exponential type $R$ and $\mathcal{F}f \in L^2(\R)$ if $f \in L^2(\R)$ and $\supp(f) \subseteq [-R, R]$.
  Conversely, if $F$ is an entire function of exponential type $R$ and $F \in L^2(\R)$, then $F = \mathcal{F} f$ for some $f \in L^2(\R)$ such that $\supp(f) \subseteq [-R, R]$.
\end{theorem}

\begin{proof}[Proof of \Cref{thm:main-theorem}]
    We begin by verifying that $\psi$ in the Beurling--Malliavin theorem can be selected such that $\int_{\R^d} \psi(x) \dif x \neq 0$.    
    Consider first the univariate case and let $\psi \in L^2(\R)$ be a function in the Beurling--Malliavin theorem such that $\supp(\psi) \subseteq [-\varepsilon, \varepsilon]$ and \smash{$\lvert \mathcal{F}\psi(\xi) \rvert \leq e^{-\omega(\xi)}$}.
    Because \smash{$\int_{\R} \psi(x) \dif x = \mathcal{F}\psi(0)$}, let us suppose to the contrary that \smash{$\mathcal{F}\psi(0) = 0$}.
    It follows from the Paley--Wiener theorem that $\mathcal{F} \psi$ is an entire function of exponential type $\varepsilon$.
    Therefore $\mathcal{F}\psi(\xi) = \sum_{n=0}^\infty a_n \xi^n$ for some $a_n \in \C$ and all $\xi \in \C$.
    Since $\mathcal{F}\psi(0) = 0$, $\mathcal{F}\psi$ has a zero of order $m \geq 1$ at the origin, meaning that $a_0 = \ldots = a_{m-1} = 0$ and $a_m \neq 0$.
    It follows that the function
    \begin{equation*}
      \phi(\xi) = \frac{\mathcal{F}\psi(\xi)}{\xi^m} = a_m + a_{m+1} \xi + a_{m+2} \xi^2 + \cdots
    \end{equation*}
    is entire and satisfies $\phi(0) \neq 0$.
    Since $\lvert \phi(\xi) \rvert \leq \lvert \mathcal{F} \psi(\xi) \rvert$ for $\lvert \xi \rvert \geq 1$ and $\phi$ is bounded on the unit disk, we conclude that $\phi$ is of exponential type $\varepsilon$, in $L^2(\R)$, and $\lvert \delta \phi(\xi) \rvert \leq e^{-\omega(\xi)}$ for some $\delta > 0$ and all $\xi \in \R$.
    Let $\varphi = \mathcal{F} (\delta \phi)$.
    Then $\varphi \in L^2(\R)$ and  $\supp(\varphi) \subseteq [-\varepsilon, \varepsilon]$ by the Paley--Wiener theorem.
    When $d = 1$, we have thus constructed a function $\varphi$ that has all the properties of $\psi$ and moreover satisfies $\int_{\R} \varphi(x) \dif x = \phi(0) \neq 0$.
    In his proof of the Beurling--Malliavin theorem for $d > 1$, Bergman~\cite{Bergman2025} constructs $\psi$ by taking a function $\psi_0$ such that $\mathcal{F}\psi_0(\xi) = \sum_{n=0}^\infty a_{2n} \xi^{2n}$ obtained from the univariate Beurling--Malliavin theorem and sets
    \begin{equation*}
      \mathcal{F} \psi(\xi) = \sum_{n=0}^\infty a_{2n} (\xi_1^2 + \cdots + \xi_d^2)^{2n}.
    \end{equation*}
    Therefore $\int_{\R^d} \psi(x) \dif x = \mathcal{F}\psi(0) = a_0 = \mathcal{F}\psi_0(0) = \int_{\R} \psi_0(x) \dif x$, so that selecting $\psi_0$ such that $\int_{\R} \psi_0(x) \dif x \neq 0$ yields $\psi$ such that $\int_{\R^d} \psi(x) \dif x \neq 0$.

    By~\eqref{eq:ass-1} any function $\psi$ from the Beurling--Malliavin theorem satisfies
    \begin{equation*}
      \int_{\R^d} \frac{\lvert \mathcal{F}\psi(\xi) \rvert^2}{\mathcal{F}\Phi(\xi)} \dif \xi \leq \int_{\R^d} \frac{e^{-2\omega(\xi)}}{\mathcal{F}\Phi(\xi)} \dif \xi \leq c^2 \int_{\R^d} \mathcal{F}\Phi(\xi) \dif \xi < \infty 
    \end{equation*}
    and is therefore an element of $H(\R^d)$ by the characterisation~\eqref{eq:rkhs}.
    Above we showed that this function can be selected such that $\int_{\R^d} \psi(x) \dif x \neq 0$.
    From $K$ being translation-invariant it follows that $\psi(\cdot - x) \in H(\R^d)$ for any $x \in \R^d$, so that for every $B(x, \varepsilon)$ there is a non-trivial $\psi \in H(\R^d)$ such that $\supp(\psi) \subseteq B(x, \varepsilon)$ and $\psi(x) = 1$.
\end{proof}

Let us then prove \Cref{thm:polynomials}.
Similar arguments will appear in~\cite{KarvonenPronzatoZhigljavsky2026}.
Part~(i) is a consequence of the following slightly more general proposition.

\begin{proposition} \label{prop:unbounded}
  Let $X$ be a normed vector space.
  Suppose that $K(x, y) = \Phi(x - y)$ for a positive-semidefinite function $\Phi \colon X \to \R$ such that $\Phi(x) \to 0$ as $\lVert x \rVert \to \infty$.
  If $\Omega \subseteq X$ is unbounded, then $f(x) \to 0$ as $\lVert x \rVert \to \infty$ for every $f \in H(\Omega)$. 
\end{proposition}
\begin{proof}
  Let $f \in H(\Omega)$ and assume to the contrary that there are $r > 0$ and a sequence $(x_n)_{n=1}^\infty$ in $\Omega$ such that $\lVert x_n \rVert \to \infty$ as $n \to \infty$ and $\lvert f(x_n) \rvert \geq r$ for all $n$.
  By passing to a subsequence and considering $-f$ if necessary we may assume that $\lVert x_{n+1} - x_n \rVert \to \infty$ as $n \to \infty$ and $f(x_n) \geq r $ for all $n$.
  Because $f \in H(\Omega)$, there is $c > 0$ such that $R(x, y) = K(x, y) - c^2 f(x) f(y) = \Phi(x - y) - c^2 f(x) f(y)$ defines a positive-semidefinite kernel on $\Omega$~\cite[Theorem~3.11]{Paulsen2016}.
  Therefore the quadratic form
  \begin{equation*}
    \begin{split}
      r_{N,l} = \sum_{n,m=1}^N R(x_{n+l}, x_{m+l}) 
      = \sum_{n,m=1}^N \big( \Phi(x_{n+l} - x_{m+l}) - c^2 f(x_{n+l}) f(x_{m+l}) \big)      
      \end{split}
  \end{equation*}
  is non-negative for every $N \geq 1$ and $l \geq 0$.
  Since $\Phi(x) \to 0$ as $\lVert x \rVert \to \infty$,
  \begin{equation*}
    \max_{\substack{ n, m \leq N \\ n \neq m}} \lvert \Phi(x_{n+l} - x_{m+l}) \rvert \leq \frac{1}{2} c^2 r^2
  \end{equation*}
  for all sufficiently large $l$.
  Therefore, for sufficiently large $l$,
  \begin{equation*}
    \begin{split}
      r_{N,l} &= \sum_{n=1}^N \Phi(x_{n+l} - x_{n+l}) + \sum_{n \neq m} \Phi(x_{n+l} - x_{m+l}) - c^2 \sum_{n,m=1}^N f(x_{n+l}) f(x_{m+l}) \\
      &\leq \Phi(0) N + \frac{1}{2} c^2 r^2 N^2  - c^2 r^2 N^2 \\
      &= \bigg( \Phi(0) - \frac{1}{2} c^2 r^2 N \bigg) N,      
      \end{split}
  \end{equation*}
  which is negative if $N > 2\Phi(0)/(c^2 r^2)$.
  It follows that $r_{N,l}$ is negative for sufficiently large $N$ and $l$ which contradicts the assumption that $f \in H(\Omega)$.
\end{proof}

\begin{proof}[Proof of \Cref{thm:polynomials}]
  Part~(i) follows directly from \Cref{prop:unbounded}.
  Consider Part~(ii).
  Suppose to the contrary that $H(\Omega)$ contains a non-trivial polynomial $p$.
  Because polynomials are real analytic and the interior of $\Omega$ is non-empty, $p$ has a unique analytic continuation onto $\R^d$.
  Since all functions in $H(\R^d)$ are real analytic by \Cref{thm:analytic}, this continuation cannot lie in $H(\R^d)$ by Part~(i).
  Therefore $p \notin H(\Omega)$.

  Consider Part~(iii).  
  Assume that $\Omega$ is bounded and that there are $c$ and $\omega$ such that~\eqref{eq:ass-1} and~\eqref{eq:ass-2} hold.
  Let $r(x) = \chi_{[-1/2,1/2]}(x)$ be a rectangular function and $r_n$ its $n$-fold convolution (i.e., $r_1 = r$).
  Then $r_n$ is a cardinal $B$-spline of order $n$, a compactly supported function that is piecewise polynomial of degree $n-1$~\cite[Chapter~4]{Chui1992}.
  Moreover,
  \begin{equation*}
    \mathcal{F}r_n(\xi) = \frac{1}{\sqrt{2\pi}} \operatorname{sinc}^n\!\bigg(\frac{\xi}{2\pi} \bigg) .
  \end{equation*}
  There is $\varepsilon$ such that $\Omega \subseteq B(0, \varepsilon)$.  
  By multiplying, scaling, and translating such functions we can construct for any multi-index $\alpha \in \N_0^d$ a function $p_\alpha$ such that $p_\alpha(x) = x^\alpha$ on $B(0, 2\varepsilon)$.
  For any function $\psi$ supported on $B(0, \varepsilon)$,
  \begin{equation*}
      P_\alpha(x) = \int_{\R^d} \psi(t) p_\alpha(x - t) \dif t = \sum_{\beta \leq \alpha} \binom{\alpha}{\beta} x^\beta  \int_{B(0, \varepsilon)} \psi(t) (-t)^{\alpha-\beta} \dif t 
  \end{equation*}
  for $x \in B(0, \varepsilon)$.
  By \Cref{thm:main-theorem}, the function $\psi$ can be selected such that $\psi \in H(\R^d)$ and $\int_{\R^d} \psi(x) \dif x = \int_{B(0, \varepsilon)} \psi(x) \dif x \neq 0$.
  Therefore 
  \begin{equation*}
    P_\alpha(x) = x^\alpha + \text{lower-order terms} 
  \end{equation*}
  on $\Omega$.
  Because $\lvert \operatorname{sinc}(x) \rvert \leq 1$, $\lvert \mathcal{F}P_\alpha \rvert \leq \lvert \mathcal{F} \psi \rvert$ and consequently $P_\alpha \in H(\R^d)$ by the characterisation~\eqref{eq:rkhs}.
  By summing appropriate scalar multiples of $P_\alpha$ for different $\alpha$ we obtain a function that coincides with any given polynomial on $\Omega$ and is an element of $H(\R^d)$.
  Therefore $H(\Omega)$ contains all polynomials.
\end{proof}

\section*{Acknowledgements}

This work was supported by the Research Council of Finland grants 359183 (``Flagship of Advanced Mathematics for Sensing, Imaging and Modelling''), and 368086 (``Inference and approximation under misspecification'').
I thank Daniel Winkle for comments and suggestions.

\bibliographystyle{abbrv}
\setlength{\bibsep}{0pt plus 0.4ex}
\bibliography{references}

@article{BeurlingMalliavin1962,
    Author = "Beurling, A. and Malliavin, P.",
    Title = "On {F}ourier transforms of measures with compact support",
    Journal = "Acta Mathematica",
    Volume = "107",
    Number = "3--4",
    Pages = "291--309",
    Year = "1962"
}

@article{KimeldorfWahba1970,
    Author = "Kimeldorf, G. S. and Wahba, G.",
    Title = "A correspondence between {B}ayesian estimation on stochastic processes and smoothing by splines",
    Journal = "The Annals of Mathematical Statistics",
    Volume = "41",
    Number = "2",
    Pages = "495--502",
    Year = "1970"
}

@book{ReedSimon1975,
    Author = "Reed, M. and Simon, B.",
    Title = "Methods of Modern Mathematical Physics. {II}: {F}ourier Analysis, Self-Adjointness",
    Publisher = "Elsevier",
    Year = "1975"
}

@book{Zolotarev1986,
    Author = "Zolotarev, V. M.",
    Title = "One-Dimensional Stable Distributions",
    Publisher = "American Mathematical Society",
    Year = "1986",
    Series = "Translations of Mathematical Monographs",
    Number = "65",
}

@book{Chui1992,
    Author = "C. K. Chui",
    Title = "An Introduction to Wavelets",
    Publisher = "Academic Press",
    Year = "1992"
}

@book{Rodino1993,
    Author = "Rodino, L.",
    Title = "Linear Partial Differential Operators in {G}evrey Spaces",
    Publisher = "World Scientific",
    Year = "1993"
}

@book{Strichartz1994,
    Author = "Strichartz, R.",
    Title = "A Guide to Distribution Theory and {F}ourier Transforms",
    Publisher = "CRC Press",
    Year = "1994"
}

@book{HavinJoricke1994,
    Author = "Havin, V. and Jöricke, B.",
    Title = "The Uncertainty Principle in Harmonic Analysis",
    Publisher = "Springer",
    Year = "1994",
    Series = " Ergebnisse der Mathematik und ihrer Grenzgebiete. 3. Folge",
    Number = "28",
}

@book{KrantzParks2002,
    Author = "Krantz, S. G. and Parks, H. R.",
    Title = "A Primer of Real Analytic Functions",
    Publisher="Birkhäuser",
    Edition = "2nd",
    Year = "2002"
}

@article{GaroniFrankel2002,
    Author = "Garoni, T. M. and Frankel, N. E.",
    Title = "{L}évy flights: Exact results and asymptotics beyond all orders Available to Purchase ",
    Journal = "Journal of Mathematical Physics",
    Volume = "43",
    Pages = "2670--2689",
    Year = "2002"
}

@book{Berlinet2004,
    Author = "Berlinet, A. and Thomas{-}Agnan, C.",
    Title = "Reproducing Kernel {H}ilbert Spaces in Probability and Statistics",
    Publisher = "Springer",
    Year = "2004"
}

@book{Wendland2005,
    Author = "Wendland, H.",
    Title = "Scattered Data Approximation",
    Publisher = "Cambridge University Press",
    Year = "2005",
    Series = "Cambridge Monographs on Applied and Computational Mathematics",
    Number = "17",
}

@Article{Mashreghi2006,
    Author = "Mashreghi, J. and Nazarov, F. L. and Havin, V. P.",
    Title = "{B}eurling--{M}alliavin multiplier theorem: {T}he seventh proof",
    Journal = "St. Petersburg Mathematical Journal",
    Volume = "17",
    Pages ="699--744",
    Year = "2006"
}

@book{RasmussenWilliams2006,
    Author = "Rasmussen, C. E. and Williams, C. K. I.",
    Title = "Gaussian Processes for Machine Learning",
    Publisher = "MIT Press",
    Series = "Adaptive Computation and Machine Learning",
    Year = "2006"
}

@article{Steinwart2006,
    Author = "Steinwart, I. and Hush, D. and Scovel, C.",
    Title = "An explicit description of the reproducing kernel {H}ilbert spaces of {G}aussian {RBF} kernels",
    Journal = "IEEE Transactions on Information Theory",
    Volume = "52",
    Number = "10",
    Pages = "4635--4643",
    Year = "2006"
}

@article{SunZhou2008,
    Author = "Sun, H.-W. and Zhou, D.-X.",
    Title = "Reproducing kernel {H}ilbert spaces associated with analytic translation-invariant {M}ercer kernels",
    Journal = "Journal of Fourier Analysis and Applications",
    Volume = "14",
    Number = "1",
    Pages = "89--101",
    Year = "2008"
}

@book{NovakWozniakowski2008,
    Author = "Novak, E. and Wo{\'z}niakowski, H.",
    Title = "Tractability of Multivariate Problems. Volume {I}: Linear Information",
    Publisher = "European Mathematical Society",
    Series = "EMS Tracts in Mathematics",
    Number = "6",
    Year = "2008"
}

@article{Minh2010,
    Author = "Minh, H. Q.",
    Title = "Some properties of {G}aussian reproducing kernel {H}ilbert spaces and their implications for function approximation and learning theory",
    Journal = "Constructive Approximation",
    Volume = "32",
    Number = "2",
    Pages = "307--338",
    Year = "2010"
}

@article{VazquezBect2010a,
    Author = "Vazquec, E. and Bect, J.",
    Title = "Convergence properties of the expected improvement algorithm with fixed mean and covariance functions",
    Journal = "Journal of Statistical Planning and Inference",
    Volume = "140",
    Number = "11",
    Pages = "3088--3095",
    Year = "2010"
}

@book{Paulsen2016,
    Author="Paulsen, V. I. and Raghupathi, M.",
    Title="An Introduction to the Theory of Reproducing Kernel {H}ilbert Spaces",
    Publisher="Cambridge University Press",
    Year = "2016"
}

@article{Mityagin2020,
    Author = "Mityagin, B. S.",
    Title = "The zero set of a real analytic function",
    Journal = "Mathematical Notes",
    Volume = "107",
    Pages = "529--530",
    Year = "2020"
}

@Article{Vasileyv2025,
  Author = "Vasilyev, I.",
  Title = "The {B}eurling and {M}alliavin theorem in several dimensions",
  Journal = "Mathematische Annalen ",
  Volume = "391",
  Pages = "6057--6072",
  Year = "2025"
}

@Article{Bergman2025,
    Author = "Bergman, A.",
    Title = "A remark on the {B}eurling-{M}alliavin theorem in several variables",
    Journal = "arXiv:2512.07271v1",
    Year = "2025"
}

@book{KarvonenPronzatoZhigljavsky2026,
    Author="Karvonen, T. and Pronzato, L. and Zhigljavsky, A.",
    Title="Space-Filling Design and Kernels: Theory and Algorithms",
    Publisher = "Springer",
    Year = "2026+"
}

\end{document}